\documentclass[envcountsame,final]{svjms}
\usepackage{graphicx}
\usepackage{times,cyrillic}
\usepackage{amssymb,amsmath2000}

\newcommand\Wr{\operatorname{Wr}}
\newcommand\Tw{\operatorname{Tw}}
\newcommand\MCN{\operatorname{Cr}}
\newcommand\Area{\operatorname{Area}}
\newcommand\Len{\operatorname{Len}}

\newcommand\Order{\mathcal{O}}
\newcommand\ds{\,\operatorname{ds}}

\newcommand\Peri[1]{P_{#1}}
\newcommand\Reals{{\mathbb R}}
\newcommand\R{\Reals}
\newcommand\Z{{\mathbb Z}}
\newcommand\Ltype{{\mathcal L}}
\newcommand\Lk{\operatorname{Lk}}
\newcommand\Flux{\operatorname{Flux}}
\newcommand\cross{\times}
\newcommand\Tube{T}
\newcommand\dArea{\,\operatorname{dA}}

\newcommand\area{\operatorname{Area}}
\newcommand\Over{\operatorname{Ov}}
\newcommand\bridge{\operatorname{Br}}

\newcommand\bdy{\partial}

\newcommand\setm{\smallsetminus}
\renewcommand\cong{\equiv}
\newcommand\half{\tfrac12}
\newcommand\ac{\operatorname{AC}}
\newcommand\pov{\operatorname{PC}}
\newcommand\genus{\operatorname{genus}}
\let\paragraph=\S
\renewcommand\S{{\mathbb S}}

\def\incgrh#1#2{\includegraphics[height=#2in]{figs/#1}}
\def\incgrs#1#2{\includegraphics[scale=#2]{figs/#1}}

\def\figside#1#2#3#4#5#6%
{ \begin{figure}\sidecaption
  \hbox {\hspace{#5in}\incgrs{#1}{#2}\hspace{#5in}}
  \caption[#3]{#4\vspace{#6in}}\label{fig:#1}
\end{figure} }

\def\figs#1#2#3#4%
{ \begin{figure}
  \begin{center} \incgrs{#1}{#2} \end{center}
  \caption[#3]{#4}\label{fig:#1}
\end{figure} }

\def\figh#1#2#3#4%
{ \begin{figure}
  \begin{center} \incgrh{#1}{#2} \end{center}
  \caption[#3]{#4}\label{fig:#1}
\end{figure} }

\def\figtwo#1#2#3#4%
{ \begin{figure}
  \begin{center} \incgrh{#1A}{#2} \hspace{.5in} \incgrh{#1B}{#2} \end{center}
  \caption[#3]{#4}\label{fig:#1}
\end{figure} }

\def\figtwox#1#2#3#4#5%
{ \begin{figure}
{\hfill\hbox to 4in
 {\vbox {\vspace{.4in}\hbox{\incgrh{#1A}{#2}} \vspace{.4in}} \hfill
  \vbox{\hbox{\incgrh{#1B}{#3}}} } \hfill}
\caption[#4]{#5}\label{fig:#1}
\end{figure} }

\def\figthree#1#2#3#4%
{ \begin{figure}
  \begin{center} \incgrh{#1A}{#2} \incgrh{#1B}{#2} \incgrh{#1C}{#2} \end{center}  \caption[#3]{#4}\label{fig:#1}
\end{figure} }

\begin{document}

\title{On the Minimum Ropelength of Knots and Links}

\author{Jason Cantarella\inst{1}
   \and Robert B. Kusner\inst{2}
   \and John M. Sullivan\inst{3}}
\institute{University of Georgia, Athens, \email{cantarel@math.uga.edu}
      \and University of Massachusetts, Amherst, \email{kusner@math.umass.edu}
      \and University of Illinois, Urbana, \email{jms@math.uiuc.edu}}
\date{Received: April 2, 2001, and in revised form February 22, 2002.}

\authorrunning{Cantarella, Kusner and Sullivan}
\maketitle

\begin{abstract}
 The ropelength of a knot is the quotient of its length by its thickness,
 the radius of the largest embedded normal tube around the knot.
 We prove existence and regularity for ropelength minimizers in
 any knot or link type; these are $C^{1,1}$ curves, but need not be
 smoother.  We improve the lower bound for the ropelength of a
 nontrivial knot, and establish new ropelength bounds for small knots
 and links, including some which are sharp.
\end{abstract}

\section*{Introduction}

How much rope does it take to tie a knot?  We measure the
\emph{ropelength} of a knot as the quotient of its length and its
\emph{thickness}, the radius of the largest embedded normal tube
around the knot.  A ropelength-minimizing configuration of
a given knot type is called \emph{tight}.

Tight configurations make interesting choices for canonical representatives
of each knot type, and are also referred to as ``ideal knots''.
It seems that geometric properties of tight knots and links
are correlated well with various physical properties of knotted polymers.
These ideas have attracted special attention in biophysics,
where they are applied to knotted loops of DNA.
Such knotted loops are important tools
for studying the behavior of various enzymes known as topoisomerases.
For information on these applications, see for instance
\cite{Sumners,Stasiak,Katritch,Stasiak-comp,DarSum1,DarSum2,CKS-nature,StasiakLaurie}
and the many contributions to the book \emph{Ideal Knots}~\cite{ideal}.

In the first section of this paper, we show the equivalence of
various definitions that have previously been given for thickness.
We use this to demonstrate that in any knot or link type
there is a ropelength minimizer, and that minimizers are necessarily
$C^{1,1}$ curves (Theorem~\ref{thm:min}).

The main results of the paper are several new lower bounds for ropelength,
proved by considering intersections of the normal tube and a spanning surface.
For a link of unit thickness, if one component is linked to~$n$
others, then its length is at least $2\pi+\Peri{n}$, where~$\Peri{n}$
is the length of the shortest curve surrounding $n$ disjoint
unit-radius disks in the plane (Theorem~\ref{thm:peribound}).
This bound is sharp in many simple
cases, allowing us to construct infinite families of tight
links, such as the simple chain shown in
Figure~\ref{fig:tight-chain}.  The only previously known example of a
tight link was the Hopf link built from two round circles, which was the
solution to the Gehring link problem~\cite{EdSc,Oss,Gage}.  Our new
examples show that ropelength minimizers need not be $C^2$, and need
not be unique.

\figh{tight-chain}{1.3}{A tight link.}
{A simple chain of~$k\ge2$ rings, the connect sum of $k-1$ Hopf links,
can be built from stadium curves (with circles at the ends).
This configuration has ropelength $(4\pi+4)k-8$
and is tight by Theorem~\ref{thm:peribound};
it shows that ropelength minimizers need not be~$C^2$.}

Next, if one component in a unit-thickness link has total linking
number~$n$ with the other components, then its length is at least
$2\pi+ 2\pi\sqrt n$, by Theorem~\ref{thm:link}.
We believe that this bound is never sharp for $n>1$.
We obtain it by using a calibration argument to estimate the area of a cone
surface spanning the given component, and the isoperimetric inequality
to convert this to a length bound.
For links with linking number zero, we need a different approach:
here we get better ropelength bounds (Theorem~\ref{thm:acbounds}) in terms of
the \emph{asymptotic crossing number} of Freedman and He~\cite{Freedman-He}.

Unit-thickness knots have similar lower bounds on length,
but the estimates are more intricate and rely on two
additional ideas.  In Theorem~\ref{thm:fourpi},
we prove the existence of a point from which
any nontrivial knot has cone angle at least $4\pi$.
In Section~\ref{sec:overx}, we
introduce the \emph{parallel overcrossing number} of a knot,
which measures how many times it crosses over a parallel
knot: we conjecture that this equals the crossing number, and we prove it
is at least the bridge number (Proposition~\ref{prop:bridge}).
Combining these ideas, we show (Theorem~\ref{thm:newbounds}) that
any nontrivial knot has ropelength at least $4\pi+ 2\pi\sqrt 2 \approx
6.83\pi \approx 21.45$.  The best previously known lower
bound~\cite{lsdr} was $5 \pi\approx 15.71$.
Computer experiments~\cite{SDKP-ideal} using Pieranski's SONO
algorithm~\cite{Pieranski}
suggest that the tight trefoil has ropelength around $32.66$. Our
improved estimate still leaves open the old question of whether any
knot has ropelength under $24$, that is: Can a knot be tied in one
foot of one-inch (diameter) rope?

\section{Definitions of Thickness}

To define the thickness of a curve, we follow the paper~\cite{GM} of
Gonzalez and Maddocks.  Although they considered only smooth curves,
their definition (unlike most earlier ones, but see \cite{KS-disto})
extends naturally to the more general curves we will need.
In fact, it is based on Menger's notion (see~\cite[\paragraph 10.1]{BluMen})
of the three-point curvature of an arbitrary metric space.

For any three distinct points $x$, $y$, $z$ in $\R^3$, we let $r(x,y,z)$ be
the radius of the (unique) circle through these points (setting
$r=\infty$ if the points are collinear).  Also, if $V_x$ is a line
through~$x$, we let $r(V_x,y)$ be the radius of the circle through $y$
tangent to $V_x$ at $x$.

Now let $L$ be a link in $\R^3$, that is, a disjoint union of simple
closed curves.  For any $x\in L$, we define the \emph{thickness} $\tau(L)$
of $L$ in terms of a \emph{local thickness} $\tau_x(L)$ at $x\in L$:
\[
  \tau(L):=\inf_{x\in L} \tau_x(L), \qquad
  \tau_x(L):=\!\!\inf_{\substack{y,z\in L\\x\ne y\ne z\ne x}} \!\!r(x,y,z).
\]
To apply this definition to nonembedded curves, note that we consider
only triples of distinct points $x,y,z\in\R^3$.  We will see later
that a nonembedded curve must have zero thickness unless its image
is an embedded curve, possibly covered multiple times.

Note that any sphere cut three times by $L$ must have radius greater
than $\tau(L)$.  This implies that the closest distance between any
two components of~$L$ is at least~$2\tau$, as follows:
Consider a sphere whose diameter achieves this minimum distance;
a slightly larger sphere is cut four times.

We usually prefer not to deal explicitly with our space curves as
maps from the circle.  But it is important to note that below,
when we talk about curves being in class $C^{k,\alpha}$, or
converging in $C^k$ to some limit, we mean with respect to the
constant-speed parametrization on the unit circle.

Our first two lemmas give equivalent definitions of thickness.
The first shows that the infimum in the
definition of~$\tau(L)$ is always attained in a limit when (at least) two
of the three points approach each other.  Thus, our definition
agrees with one given earlier by Litherland \emph{et~al.}~\cite{lsdr}
for smooth curves.  If $x\ne y\in L$ and $x-y$ is perpendicular
to both $T_xL$ and $T_yL$, then we call $|x-y|$ a \emph{doubly critical
self-distance} for $L$.

\begin{lemma}\label{lem:gm-meets-lsdr}
  Suppose $L$ is $C^1$, and let $T_x L$ be its tangent line at $x\in L$. Then 
the thickness is given by
\[
    \tau(L) = \!\inf_{x\ne y\in L} r(T_x L,y).
\]
  This equals the infimal radius of curvature of~$L$ or half
  the infimal doubly critical self-distance, whichever is less.
\end{lemma}
\begin{proof}
  The infimum in the definition of thickness either is achieved for
  some distinct points $x$, $y$, $z$, or is approached along the
  diagonal when $z$, say, approaches $x$, giving us $r(T_xL,y)$.  But the first
  case cannot happen unless the second does as well: consider the
  sphere of radius $\tau=r(x,y,z)$ with $x$, $y$ and $z$ on its
  equator, and relabel the points if necessary so that $y$ and $z$ are
  not antipodal.  Since this $r$ is infimal, $L$ must be tangent to
  the sphere at $x$.  Thus $r(T_xL,y)\le\tau$, and we see that $\tau(L)
  = \inf r(T_x L,y)$.  This infimum, in turn, is achieved either for
  some $x\ne y$, or in a limit as $y\to x$ (when it is the infimal
  radius of curvature).  In the first case, we can check that $x$ and
  $y$ must be antipodal points on a sphere of radius $\tau$, with $L$
  tangent to the sphere at both points.  That means, by definition,
  that $2\tau$ is a doubly critical self-distance for $L$.
\qed\end{proof}

A version of Lemma~\ref{lem:gm-meets-lsdr} for smooth curves appeared
in~\cite{GM}.  Similar arguments there show that the local thickness
can be computed as
\[
\tau_x(L)=\inf_{y\ne x} r(T_y L,x).
\]

\begin{lemma}\label{lem:tau=reach=nir}
  For any $C^1$ link $L$, the thickness of~$L$ equals the reach of~$L$;
  this is also the normal injectivity radius of~$L$.
\end{lemma}

The \emph{reach} of a set $L$ in $\R^3$, as defined by
Federer~\cite{Federer}, is the largest~$\rho$ for which any point~$p$ in
the $\rho$-neighborhood of~$L$ has a unique nearest point in $L$.  The
\emph{normal injectivity radius} of a $C^1$ link $L$ in $\R^3$
is the largest~$\iota$ for which the union of the open normal
disks to~$L$ of radius $\iota$ forms an embedded tube.

\begin{proof}
  Let $\tau$, $\rho$, and $\iota$ be the thickness, reach, and normal
  injectivity radius of~$L$.  We will show that
  $\tau\leq \rho\leq \iota\leq\tau$.

  Suppose some point $p$ has two nearest neighbors
  $x$ and $y$ at distance~$\rho$.  Thus $L$ is tangent at $x$ and $y$ to
  the sphere around~$p$, so a nearby sphere cuts~$L$ four
  times, giving $\tau \leq \rho$.

  Similarly, suppose some~$p$ is on two normal circles of~$L$ of radius
  $\iota$.  This~$p$ has two neighbors on~$L$ at distance~$\iota$,
  so $\rho\leq \iota$.

  We know that $\iota$ is less than the infimal radius of curvature of
  $L$. Furthermore, the midpoint of a chord of~$L$ realizing the infimal
  doubly self-critical distance of~$L$ is on two normal disks of~$L$.
  Using Lemma~\ref{lem:gm-meets-lsdr}, this shows that $\iota \leq \tau$,
  completing the proof.
\qed\end{proof}

If $L$ has thickness $\tau>0$, we will call the embedded (open) normal tube
of radius $\tau$ around~$L$ the \emph{thick tube} around~$L$.

We define the \emph{ropelength} of a link $L$ to be $\Len(L)/\tau(L)$, the
(scale-invariant) quotient of length over thickness.  Every curve
of finite roplength is $C^{1,1}$, by Lemma~\ref{lem:sec-is-lip} below.
Thus, we are free to restrict our attention to~$C^{1,1}$ curves,
rescaled to have (at least) unit thickness.  This means they have
embedded unit-radius normal tubes, and curvature bounded above
by~$1$.  The ropelength of such a curve is (at most) its length.

\section{Existence and Regularity of Ropelength Minimizers}\label{sec:exist}

We want to prove that, within every knot or link type, there exist
curves of minimum ropelength.  The lemma below allows us to use the
direct method to get minimizers.  If we wanted to, we could work with
$C^1$ convergence in the space of $C^{1,1}$ curves, but it seems better to
state the lemma in this stronger form, applying to all rectifiable
links.

\begin{lemma}\label{lem:thi-is-usc}
  Thickness is upper semicontinuous with respect to the $C^0$ topology
  on the space of $C^{0,1}$ curves.
\end{lemma}
\begin{proof}
  This follows immediately from the definition, since $r(x,y,z)$ is a
  continuous function (from the set of triples of distinct points in space)
  to~$(0,\infty]$.  For, if curves $L_i$ approach $L$,
  and $r(x,y,z)$ nearly realizes the thickness of~$L$, then nearby
  triples of distinct points bound from above the thicknesses of the $L_i$.
\qed\end{proof}

This proof (compare~\cite{KS-disto}) is essentially the same as the
standard one for the lower semicontinuity of length, when length of
an arbitrary curve is defined as the
supremal length of inscribed polygons.  Note that thickness can jump
upwards in a limit, even when the convergence is $C^1$.  For instance,
we might have an elbow consisting of two straight segments connected
by a unit-radius circular arc whose angle decreases to zero, as shown
in Figure~\ref{fig:thi-is-semicontinuous}.

\figs{thi-is-semicontinuous}{1}{A sequence of curves with unit
thickness converge to a curve with infinite thickness} {These elbow
curves consist of straight segments connected by a circular arc of
unit radius.  No matter how small the angle of the arc is, the curve
has unit thickness, as demonstrated by the maximal embedded normal
neighborhoods shown.  These $C^{1,1}$ curves converge in $C^{1}$ to a
straight segment, with infinite thickness.}

When minimizing ropelength within a link type, we care only about
links of positive thickness $\tau>0$.  We next prove three lemmas about
such links.  It will be useful to consider the \emph{secant map}
$S$ for a link $L$, defined, for $x\ne y\in L$, by
\[ 
  S(x,y) := \pm \frac{x - y}{|x - y|} \in \R P^2.
\] 
Note that as $x\to y$, the limit of $S(x,y)$, if it exists,
is the tangent line $T_yL$.  Therefore,
the link is $C^1$ exactly if $S$ extends continuously
to the diagonal~$\Delta$ in~$L\cross L$, and is~$C^{1,1}$ exactly when this
extension is Lipschitz.  When speaking
of particular Lipschitz constants we use the following metrics:
on~$L\cross L$, we sum the (shorter) arclength distances in the
factors; on~$\R P^2$ the distance between two points is
$d=\sin\theta$, where $\theta$ is the angle between (any) lifts
of the points to~$\S^2$.

\begin{lemma}\label{lem:sec-is-lip}
  If~$L$ has thickness $\tau>0$, then its secant map~$S$
  has Lipschitz constant $1/2\tau$.  Thus $L$ is~$C^{1,1}$.
\end{lemma}
\begin{proof}
  We must prove that $S$ has Lipschitz constant $1/2\tau$ on
  $(L\cross L)\setm\Delta$; it then has a Lipschitz extension.
  By the triangle inequality, it suffices to
  prove, for any fixed $x \in L$, that
  $d\big(S(x,y),S(x,z)\big) \leq |y-z|/2\tau$
  whenever $y$ and $z$ are sufficiently close along $L$.
  Setting $\theta:=\angle zxy$, we have
  \[ 
    d\big(S(x,y),S(x,z)\big) = \sin\theta
     = \frac {|y-z|}{2\,r(x,y,z)} \leq \frac {|y-z|}{2\tau},
  \]   
  using the law of sines and the definition $\tau := \inf r$.
\qed\end{proof}
\figside{secant-map}{1.2}{The secant map.}
 {The secant map for a thick knot is Lipschitz by Lemma~\ref{lem:sec-is-lip}:
  when $y$ and $z$ are close along the knot, the secant directions
  $\overline{xy}$ and $\overline{xz}$ are close.  Here $r=r(x,y,z)$
  is an upper bound for the thickness of the knot.}{0}{.1}

Although we are primarily interested in links (embedded curves), we note that
Lemma~\ref{lem:sec-is-lip} also shows that a nonembedded curve $L$ must have
thickness zero, unless its image is contained in some embedded curve.
For such a curve $L$ contains some point $p$ where at least three arcs
meet, and at least one pair of those arcs will fail to join in a $C^{1,1}$
fashion at~$p$.

\begin{lemma}\label{lem:secants}
  If $L$ is a link of thickness $\tau > 0$, then any points $x,y\in L$
  with $|x-y|<2\tau$ are connected by an arc of~$L$ of length at most
 $$2\tau\arcsin\frac{|x-y|}{2\tau} \,\leq \,\frac\pi2 \,|x-y|.$$
\end{lemma}
\begin{proof}
  The two points $x$ and $y$ must be on the same component of~$L$,
  and one of the arcs of~$L$ connecting them is contained in the
  ball with diameter~$\overline{xy}$.  By
  Lemma~\ref{lem:gm-meets-lsdr}, the curvature of~$L$ is less than
  $1/\tau$.  Thus by Schur's lemma,
  the length of this arc of~$L$ is at most $2\tau \arcsin (|x-y|/2\tau)$,
  as claimed.  Note that Chern's proof~\cite{Chern} of Schur's lemma
  for space curves, while stated only for $C^2$ curves, applies directly
  to $C^{1,1}$ curves, which have Lipschitz tantrices on the unit sphere.
  (As Chern notes, the lemma actually applies even to curves with corners,
  when correctly interpreted.)
\qed\end{proof}

\begin{lemma}
\label{lem:c0-implies-c1}
Suppose $L_i$ is a sequence of links of thickness at least $\tau>0$,
converging in $C^0$ to a limit link~$L$.  Then the convergence is
actually $C^1$, and $L$ is isotopic to (all but finitely many of)
the~$L_i$.
\end{lemma}
\begin{proof}
  To show $C^1$ convergence, we will show that the secant maps of the~$L_i$
  converge (in $C^0$) to the secant map of~$L$.  Note that when
  we talk about convergence of the secant maps, we view them (in terms
  of constant-speed parametrizations of the $L_i$) as maps from a common
  domain.  Since these maps are uniformly Lipschitz,
  it suffices to prove pointwise convergence.

  So consider a pair of points $p$, $q$ in~$L$.  Take $\epsilon < |p-q|$.
  For large enough~$i$, $L_i$~is within
  $\epsilon^2/2$ of $L$ in $C^0$, and hence the corresponding points
  $p_i$, $q_i$ in~$L_i$ have $|p_i - p| < \epsilon^2/2$ and $|q_i - q| <
  \epsilon^2/2$.
  We have moved the endpoints of the segment $\overline{pq}$
  by relatively small amounts,
  and expect its direction to change very little. In fact, the angle
  $\theta$ between $p_i - q_i$ and $p - q$ satisfies
  $\sin \theta < (\epsilon^2/2 + \epsilon^2/2)/\epsilon = \epsilon$.
  That is, the distance in $\R P^2$ between
  the points $S_i(p_i,q_i)$ and $S(p,q)$ is given by $\sin\theta<\epsilon$.

  Therefore, the secant maps converge pointwise, which shows 
  that the~$L_i$ converge in~$C^1$ to~$L$.  Since the
  limit link $L$ has thickness at least~$\tau$ by
  Lemma~\ref{lem:thi-is-usc}, it is surrounded by an embedded normal
  tube of diameter~$\tau$. Furthermore, all (but finitely many) of the~$L_i$
  lie within this tube, and by $C^1$ convergence are transverse
  to each normal disk. Each such $L_i$ is isotopic to~$L$ by a
  straight-line homotopy within each normal disk.
\qed\end{proof}

Our first theorem establishes the existence of \emph{tight}
configurations (ropelength minimizers) for any link type.  
This problem is interesting only for tame links: a wild link has no
$C^{1,1}$ realization, so its ropelength is always infinite.

\begin{theorem}\label{thm:min}
  There is a ropelength minimizer in any (tame) link type;
  any minimizer is $C^{1,1}$, with bounded curvature.
\end{theorem}
\begin{proof}
  Consider the compact space of all $C^{1,1}$ curves of length at
  most~$1$.  Among those isotopic to a given link~$L_0$,
  find a sequence $L_i$ supremizing the thickness.  The lengths
  of~$L_i$ approach~$1$, since otherwise rescaling would give thicker curves.
  Also, the thicknesses approach some $\tau>0$, the
  reciprocal of the infimal ropelength for the link type.  Replace the
  sequence by a subsequence converging in the $C^1$ norm to some link $L$.
  Because length is lower semicontinuous, and thickness is upper
  semicontinuous (by Lemma~\ref{lem:thi-is-usc}), the ropelength
  of~$L$ is at most $1/\tau$. By Lemma~\ref{lem:c0-implies-c1}, all
  but finitely many of the $L_i$ are isotopic to~$L$, so $L$ is
  isotopic to~$L_0$.

  By Lemma~\ref{lem:sec-is-lip}, tight links must
  be $C^{1,1}$, since they have positive thickness.
\qed\end{proof}

This theorem has been extended by Gonzalez \emph{et~al.}~\cite{gmsvdm},
who minimize a broad class of energy functionals subject to the
constraint of fixed thickness.  See also~\cite{GdlL}.

Below, we will give some examples of tight links which show
that $C^2$ regularity cannot be expected in general,
and that minimizers need not be unique.

\section{The Ropelength of Links}\label{sec:peribound}

Suppose in a link $L$ of unit thickness, some component $K$ is
topologically linked to $n$ other components~$K_i$.  We will give a
sharp lower bound on the length of~$K$ in terms of~$n$.  When every
component is linked to $n\leq 5$ others, this sharp bound lets us construct
tight links.

To motivate the discussion below, suppose $K$ was a
planar curve, bounding some region~$R$ in the plane. Each $K_i$ would then
have to puncture $R$. Since each $K_i$ is surrounded by a unit-radius
tube, these punctures would be surrounded by disjoint disks of unit
radius, and these disks would have to avoid a unit-width ribbon around
$K$. It would then be easy to show that the length of~$K$ was
at least~$2\pi$ more than $\Peri{n}$,
the length of the shortest curve surrounding
$n$ disjoint unit-radius disks in the plane.

To extend these ideas to nonplanar curves, we need to consider cones.
Given a space curve $K$ and a point $p\in\R^3$, the \emph{cone} over~$K$
from~$p$ is the disk consisting of all line segments from $p$ to points in~$K$.
The cone is intrinsically flat away from the single cone point~$p$,
and the \emph{cone angle} is defined to be the angle obtained at~$p$ if we
cut the cone along any one segment and develop it into the Euclidean plane.
Equivalently, the cone angle is the length of the projection of~$K$ to
the unit sphere around~$p$.  Note that the total Gauss curvature
of the cone surface equals $2\pi$ minus this cone angle.

Our key observation is that every space curve may be coned to some
point $p$ in such a way that the intrinsic geometry of the cone
surface is Euclidean.  We can then apply the argument above in the
intrinsic geometry of the cone.  In fact, we can get even better
results when the cone angle is greater than $2\pi$.  We first prove a
technical lemma needed for this improvement.
Note that the lemma would remain true without
the assumptions that~$K$ is $C^{1,1}$ and has curvature at
most~$1$.  But we make use only of this case, and the more general case
would require a somewhat more complicated proof.

\begin{lemma}\label{lem:pushoff}
  Let $S$ be an infinite cone surface with cone angle $\theta\ge2\pi$
  (so that $S$ has nonpositive curvature and is intrinsically
  Euclidean away from the single cone point). Let
  $R$ be a subset of $S$ which includes the cone point, and let
  $\ell$ be a lower bound for the length of any curve in $S$ surrounding $R$.
  Consider a $C^{1,1}$ curve $K$ in $S$ with geodesic curvature
  bounded above by~$1$.  If $K$ surrounds $R$ while remaining at least
  unit distance from $R$, then $K$ has length at least $\ell + \theta$.
\end{lemma}
\begin{proof}
  We may assume that $K$ has nonnegative geodesic curvature almost everywhere.
  If not, we simply replace it by the boundary of its convex hull within $S$,
  which is well-defined since $S$ has nonpositive curvature.  This boundary
  still surrounds $R$ at unit distance, is $C^{1,1}$, and has
  nonnegative geodesic curvature.

  For $t<1$, let $K_t$ denote the inward normal pushoff,
  or parallel curve to $K$, at distance $t$ within the cone.
  Since the geodesic curvature of~$K$ is bounded by~$1$, these
  are all smooth curves, surrounding $R$ and hence surrounding the cone point.
  If~$\kappa_g$ denotes the geodesic curvature of $K_t$ in $S$,
  the formula for first variation of length is
  \[ \frac{d}{dt}\Len(K_t) = -\int_{\!K_t}\!\! \kappa_g \ds = -\theta, \]
  where the last equality comes from Gauss--Bonnet,
  since $S$ is intrinsically flat except at the cone point.
  Thus $\Len(K) = \Len(K_t) + t\theta$;
  since $K_t$ surrounds $R$ for every $t < 1$, it has length at least $\ell$,
  and we conclude that $\Len(K) \geq \ell + \theta$.
\qed\end{proof} 		

\begin{lemma}\label{lem:twopi}
  For any closed curve $K$, there is a point $p$ such that the
  cone over~$K$ from~$p$ has cone angle $2\pi$.  When $K$ has
  positive thickness, we can choose $p$ to
  lie outside the thick tube around~$K$.
\end{lemma}
\begin{proof}
  Recall that the cone angle at~$q$ is given by the length of the radial
  projection of~$K$ onto the unit sphere centered at~$q$. If we choose
  $q$ on a chord of~$K$, this projection joins two antipodal points,
  and thus must have length at least~$2\pi$.  On any doubly critical chord
  (for instance, the longest chord) the point $q$ at distance $\tau(K)$
  from either endpoint must lie outside the thick tube, by
  Lemma~\ref{lem:tau=reach=nir}.

  Note that the cone angle approaches $0$ at points far from $K$.
  The cone angle is a continuous function on the complement of $K$ in $\R^3$,
  a connected set.  When $K$ has positive thickness, even the complement
  of its thick tube is connected.
  Thus if the cone angle at $q$ is greater than $2\pi$,
  the intermediate value theorem lets us choose some $p$ (outside
  the tube) from which the cone angle is exactly $2\pi$.
  Figure~\ref{fig:3-1-2pi} shows such a cone on a trefoil knot.
\qed\end{proof}

\figtwo{3-1-2pi}{1.6}{A $2\pi$ cone on a trefoil knot.}
{Two views of the same cone, whose cone angle is precisely $2\pi$,
on a symmetric trefoil knot.}

Our first ropelength bound will be in terms of a quantity we call
$\Peri{n}$, defined to equal the shortest length of any plane curve
enclosing~$n$ disjoint unit disks.  Considering the centers of the
disks, using Lemma~\ref{lem:pushoff}, and scaling by a factor of~$2$,
we see that $\Peri{n}=2\pi+2Q_n$, where $Q_n$ is the length of the
shortest curve enclosing~$n$ points separated by unit distance in the plane.

\begin{figure}
\centering
\vbox{
\vspace{.1in}
 \begin{tabular}{r|ccc}
  {} & \incgrs{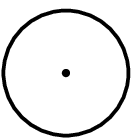}{.9} & \incgrs{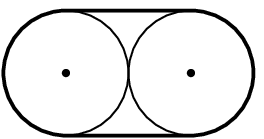}{.9}
        & \incgrs{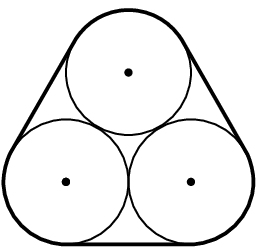}{.9} \\ \hline
  $n$        & $1$     & $2$      & $3$      \\
  $\Peri{n}$ & $2\pi$  & $2\pi+4$ & $2\pi+6$ \\
 \end{tabular}
\vskip.2in
 \begin{tabular}{r|cl}
  {} & \incgrs{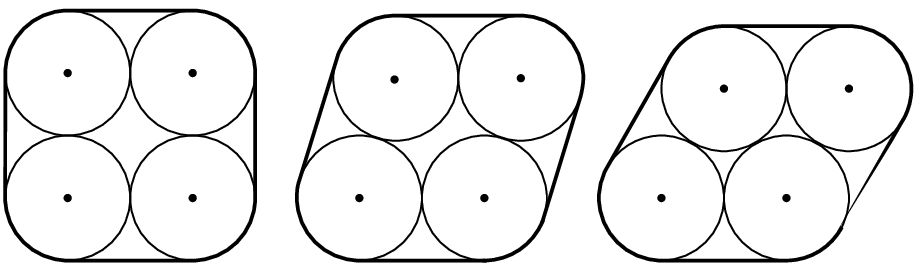}{.9}\hspace{-1.6in} & \hspace{1.3in} \\ \hline
  $n$        & $4$      &  \\
  $\Peri{n}$ & $2\pi+8$ &  \\
 \end{tabular}
\vspace{.1in}
}
\caption[Values of~$\Peri{n}$ for small $n$.]
     {The shortest curve enclosing $n$ unit disks in the plane has
      length $\Peri{n}$, and is unique for $n<4$. For $n=4$, there
      is a one-parameter family of equally short curves.}
\label{fig:peritable}
\end{figure}

For small $n$ it is not hard to determine $Q_n$ and $\Peri{n}$ explicitly
from the minimizing configurations shown in Figure~\ref{fig:peritable}.
Clearly $\Peri{1}=2\pi$, while for $2\le n\le 5$, we have
$\Peri{n}=2\pi+2n$ since $Q_n=n$.
Note that the least-perimeter curves in Figure~\ref{fig:peritable}
are unique for $n<4$, but for $n=4$ there is a continuous family of minimizers.
For $n=5$ there is a two-parameter family,
while for $n=6$ the perimeter-minimizer is again unique, with $Q_6=4+\sqrt3$.
It is clear that $Q_n$ grows like~$\sqrt n$ for $n$ large.\footnote{This
perimeter problem does not seem to have been considered previously.
However, Sch\"urmann~\cite{Sch} has also recently examined
this question.  In particular, he conjectures that the
minimum perimeter is achieved (perhaps not uniquely)
by a subset of the hexagonal circle packing
for $n<54$, but proves that this is not the case for $n>370$.}

\begin{theorem}\label{thm:peribound}
  Suppose $K$ is one component of a link of unit thickness, and the
  other components can be partitioned into $n$ sublinks, each of which
  is topologically linked to $K$. Then the length of~$K$ is at least
  $2\pi+\Peri{n}$, where $\Peri{n}$ is the minimum length of any curve
  surrounding $n$ disjoint unit disks in the plane.
\end{theorem}
\begin{proof}
  By Lemma~\ref{lem:twopi} we can find a point $p$ in space, outside the
  unit-radius tube surrounding $K$, so that coning $K$ to $p$
  gives a cone of cone angle $2\pi$, which is intrinsically flat.

  Each of the sublinks $L_i$ nontrivially linked to $K$ must puncture
  this spanning cone in some point $p_i$. Furthermore, the fact that the link
  has unit thickness implies that the $p_i$ are separated from each
  other and from~$K$ by distance at least~$2$ in space, and thus by
  distance at least~$2$ within the cone.

  Thus in the intrinsic geometry of the cone, the $p_i$ are
  surrounded by disjoint unit-radius disks, and $K$ surrounds
  these disks while remaining at least unit distance from them. Since
  $K$ has unit thickness, it is~$C^{1,1}$ with curvature bounded
  above by~$1$. Since the geodesic curvature of~$K$ on the cone
  surface is bounded above by the curvature of~$K$ in space, we can
  apply Lemma~\ref{lem:pushoff} to complete the proof.
\qed\end{proof}

For $n \leq 5$, it is easy to construct links which achieve these
lower bounds and thus must be tight.  We just ensure that each
component linking $n$ others is a planar curve of length equal to our
lower bound $2\pi+\Peri{n}$.  In particular, it must be the outer
boundary of the unit neighborhood of some curve achieving $\Peri{n}$.
In this way we construct the tight chain of
Figure~\ref{fig:tight-chain}, as well as infinite families of more
complicated configurations, including the link in
Figure~\ref{fig:complex}. These examples may help to calibrate the
various numerical methods that have been used to compute ropelength
minimizers~\cite{Pieranski,Rawdon,Laurie}. For $n\ge 6$, this
construction does not work, as we are unable to simultaneously
minimize the length of $K$ and the length of all the components it links.

\figh{complex}{1.5}{A more complicated tight link.}
{Each component in this link uses one of the perimeter-minimizing shapes
from Figure~\ref{fig:peritable}, according to how many other components it
links.  The link is therefore tight by Theorem~\ref{thm:peribound}.
This minimizer is not unique, in that some components could be rotated
relative to others.  In other examples, even the shape of individual
components (linking four others) fails to be unique.}

These explicit examples of tight links answer some existing questions
about ropelength minimizers. First, these minimizers fail, in a strong
sense, to be unique: there is a one-parameter family of tight
five-component links based on the family of curves with length
$\Peri{4}$. So we cannot hope to add uniqueness to the conclusions of
Theorem~\ref{thm:min}. In addition, these minimizers (except for the Hopf link)
are not $C^2$. This tells us that there can be no better global regularity
result than that of Theorem~\ref{thm:min}.  However, we could still hope 
that every tight link is \emph{piecewise} smooth, or even piecewise analytic.

Finally, note that the ropelength of a composite link should be somewhat
less than the sum of the lengths of its factors.  It was observed
in~\cite{Stasiak} that this deficit seems to be at least $4\pi-4$.
Many of our provably tight examples, like the simple chain in
Figure~\ref{fig:tight-chain} or the link in Figure~\ref{fig:complex},
are connect sums which give precise confirmation of this observation.

\section{Linking Number Bounds}\label{sec:linking}

We now adapt the cone surface arguments to find a lower bound on
ropelength in terms of the linking number. These bounds are more
sensitive to the topology of the link, but are not sharp, and thus
provide less geometric information. In Section~\ref{sec:acbounds}, we
will present a more sophisticated argument, which implies
Theorem~\ref{thm:link} as a consequence. However, the argument here
is concrete enough that it provides a nice introduction to the methods
used in the rest of the paper.

\begin{theorem}\label{thm:link}
  Suppose $L$ is a link of unit thickness.  If $K$ is one component
  of the link, and $J$ is the union of any collection of other components,
  let $\Lk(J,K)$ denote the total linking number of~$J$ and $K$, for
  some choice of orientations.  Then
  \[
    \Len(K) \geq 2\pi + 2\pi\sqrt{\Lk(J,K)}.
  \]
\end{theorem}
\begin{proof}
  As in the proof of Theorem~\ref{thm:peribound}, we apply
  Lemma~\ref{lem:twopi} to show that we can find
  an intrinsically flat cone surface $S$ bounded by $K$.
  We know that $K$ is surrounded by an embedded unit-radius tube $T$;
  let~$R = S \setm T$ be the portion of the cone surface outside the tube.
  Each component of~$J$ is also surrounded by an embedded unit-radius tube
  disjoint from $T$.  
  Let $V$ be the $C^1$ unit vectorfield normal to the normal disks
  of these tubes. A simple computation shows that $V$ is a
  divergence-free field, tangent to the boundary of each tube, with
  flux $\pi$ over each spanning surface inside each tube. A cohomology
  computation (compare~\cite{mutual-helicity})
  shows that the total flux of~$V$ through $R$ is
  $\Flux_R(V) = \pi \Lk(J,K)$.  Since $V$ is a unit
  vectorfield, this implies that 
  \[ 
   \Flux_R(V) = \int_R V \cdot \vec{n} \dArea
             \leq \int_R\! \dArea = \Area(R).
  \] 
  Thus $ \Area(R) \geq \pi\Lk(J,K)$.
  The isoperimetric inequality within~$S$ implies that
  any curve on $S$ surrounding $R$ has length at least
  $2\pi\sqrt{\Lk(J,K)}$. Since $L$ has unit thickness, the
  hypotheses of Lemma~\ref{lem:pushoff} are fulfilled, and we conclude that
  \[ 
   \Len(K) \geq 2\pi + 2\pi\sqrt{\Lk(J,K)},
  \] 
  completing the proof.
\qed\end{proof}

Note that the term $2\pi\sqrt{\Lk(J,K)}$ is the perimeter
of the disk with the same area as $n:=\Lk(J,K)$ unit disks.
We might hope to replace this term by $\Peri{n}$, but this seems difficult:
although our assumptions imply that $J$ punctures
the cone surface $n$ times, it is possible that there are many
more punctures, and it is not clear how to show that an appropriate
set of $n$ are surrounded by disjoint unit disks.

For a link of two components with linking number $2$,
like the one in Figure~\ref{fig:linking-number-two},
this bound provides an improvement on Theorem~\ref{thm:peribound},
raising the lower bound on the ropelength of each component to
$2\pi+2\pi\sqrt{2}$, somewhat greater than $4\pi$.

\figside{linking-number-two}{1}
{A link with ropelength at least $4\pi(1+\sqrt2)$.}
{The two components of this $(2,4)$-torus link have linking number two,
so by Theorem~\ref{thm:link}, the total ropelength is at least
$4\pi(1+\sqrt2) \approx 33.34$.  Laurie \emph{et~al.}~\cite{StasiakLaurie}
have computed a configuration with ropelength approximately $41.2$.}{0.5}{.2}

We note that a similar argument bounds the ropelength of any
curve $K$ of unit thickness, in terms of its writhe.
We again consider the flux of $V$ through a flat cone $S$.
If we perturb $K$ slightly to have rational writhe (as below in the proof of
Theorem~\ref{thm:acbounds}) and use the result that
``link equals twist plus writhe''~\cite{cal1,white},
we find that this flux is at least $|\!\Wr(K)|$, so that
\[
        \Len(K) \geq 2 \pi \sqrt{|\!\Wr(K)|}\,.
\]
There is no guarantee that this flux occurs away from the boundary of the cone,
however, so Lemma~\ref{lem:pushoff} does not apply.
Unfortunately, this bound is weaker than the
corresponding result of Buck and Simon~\cite{BuckSimon},
\[
	\Len(K) \geq 4 \pi \sqrt{|\!\Wr(K)|}\,.
\]

\section{Overcrossing Number}\label{sec:overx}

In Section~\ref{sec:linking}, we found bounds on the ropelength of links;
to do so, we bounded the area of that portion of the cone surface outside
the tube around a given component $K$ in terms of the flux of a
certain vectorfield across that portion of the surface. This argument
depended in an essential way on linking number being a signed
intersection number. 

For knots, we again want a lower bound for the area of that portion of
the cone that is at least unit distance from the boundary. But this is
more delicate and requires a more robust topological invariant.  Here,
our ideas have paralleled those of Freedman and He (see
\cite{Freedman-He,He}) in many important respects, and we adopt some
of their terminology and notation below.

Let $L$ be an (oriented) link partitioned into two parts $A$ and $B$.
The linking number $\Lk(A,B)$ is the sum of the signs of the crossings
of~$A$ over~$B$; this is the same for any projection of any link
isotopic to~$L$.  By contrast, the \emph{overcrossing number}
$\Over(A,B)$ is the (unsigned) number of crossings of~$A$ over $B$,
minimized over all projections of links isotopic to~$L$.

\begin{lemma}
\label{lem:properties-of-Ov}
For any link partitioned into two parts $A$ and~$B$,
the quantities $\Lk(A,B)$ and $\Over(A,B)$ are symmetric in~$A$ and~$B$,
and we have
$$|\Lk(A,B)| \leq \Over(A,B)\,; \qquad \Lk(A,B) \cong \Over(A,B) \pmod2\,.$$
\end{lemma}
\begin{proof} 
To prove the symmetry assertions, take any planar projection
with $n$ crossings of~$A$ over $B$. Turning the plane over,
we get a projection with $n$ crossings of~$B$ over $A$;
the signs of the crossings are unchanged.
The last two statements are immediate from
the definitions in terms of signed and unsigned sums.
\qed\end{proof}

Given a link $L$, we define its \emph{parallel overcrossing number}
$\pov(L)$ to be the minimum of $\Over(L,L')$
taken over all parallel copies $L'$ of the link $L$.
That means $L'$ must be an isotopic link such that corresponding
components of $L$ and $L'$ cobound annuli, the entire collection of
which is embedded in $\R^3$.
This invariant may be compared to Freedman and He's
\emph{asymptotic crossing number} $\ac(L)$ of~$L$, defined by
\[
\ac(L) = \inf_{pL, qL} \frac{\Over(pL,qL)}{|pq|},
\]
where the infimum is taken over all degree-$p$ satellites $pL$
and degree-$q$ satellites $qL$ of~$L$.
(This means that $pL$ lies in a solid torus around $J$ and represents
$p$~times the generator of the first homology group of that torus.)
Clearly,
\[
	\ac(L) \leq \pov(L) \leq \MCN(L),
\]
where $\MCN(L)$ is the crossing number of~$L$. 
It is conjectured that the asymptotic crossing number of~$L$ is equal
to the crossing number. This would imply our weaker conjecture:
\begin{conjecture}
If $L$ is any knot or link, $\pov(L)=\MCN(L)$.
\end{conjecture}
To see why this conjecture is reasonable, suppose $K$ is an
alternating knot of crossing number $k$.  It is
known~\cite{lt,thistle}, using the Jones polynomial, that the crossing
number of $K\cup K'$ is least $4k$ for any parallel $K'$.
It is tempting to assume that
within these $4k$ crossings of the two-component link, we can find not
only $k$ self-crossings of each knot $K$ and~$K'$, but also $k$
crossings of~$K$ over $K'$ and $k$ crossings of $K'$ over $K$.
Certainly this is the case in the standard picture of~$K$ and a planar
parallel $K'$.

Freedman and He have shown~\cite{Freedman-He} that for any knot,
\[\ac(K) \geq 2 \genus(K) - 1,\]
and hence that we have $\ac(K) \geq 1$ if $K$
is nontrivial.  For the parallel overcrossing number, our stronger
hypotheses on the topology of~$L$ and~$L'$ allow us to find a better
estimate in terms of the \emph{reduced bridge number} $\bridge(L)$.
This is the minimum number of local maxima of any height function
(taken over all links isotopic to $L$) minus the number of unknotted
split components in $L$.

\begin{proposition}
\label{prop:bridge}
For any link $L$, we have $\pov(L) \ge \bridge(L)$.
In particular, if~$L$ is nontrivial, $\pov(L)\geq 2$.
\end{proposition}
\begin{proof}
  By the definition of parallel overcrossing number, we can isotope
  $L$ and its parallel $L'$ so that, except for $\pov(L,L')$ simple
  clasps, $L'$ lies above, and $L$ lies below, a slab in $\R^3$.
  Next, we can use the embedded annuli which cobound corresponding
  components of~$L$ and $L'$ to isotope the part of~$L$ below the slab
  to the lower boundary plane of the slab.  This gives a presentation
  of~$L$ with $\pov(L)$ bridges, as in Figure~\ref{fig:clasps-and-bridges}.
\qed\end{proof}
\figh{clasps-and-bridges}{1.1}{Constructing a presentation with $\pov(L)$
  bridges.}  
{We show three stages of the proof of Proposition~\ref{prop:bridge}:
At the left, we show a projection of $L\cup L'$ with $\pov(L)$ overcrossings.
In the center, we lift $L'$ until it and $L$ lie respectively above and below
a slab, except for $\pov(L)$ simple clasps.
At the right, we isotope $L$ to flatten the undercrossings onto the boundary
of the slab and thus show that the clasps are the only bridges in $L$.}

\section{Finding a Point with Larger Cone Angle}\label{sec:coneang}

The bounds in Theorems~\ref{thm:peribound} and~\ref{thm:link} depended
on Lemma~\ref{lem:twopi} to construct a cone with cone angle $\theta = 2\pi$,
and on Lemma~\ref{lem:pushoff} to increase the total ropelength by
at least $\theta$.  For single unknotted curves, this portion of our
argument is sharp: a convex plane curve has maximum cone angle $2\pi$,
at points in its convex hull.

However, for nontrivial knots and links, we can
improve our results by finding points with greater cone angle.  In
fact, we show every nontrivial knot or link has a $4\pi$ cone point.  
The next lemma is due to Gromov~\cite[Thm.~8.2.A]{GromovFRM} and also appears
as~\cite[Thm.~1.3]{EWW}:

\begin{lemma} \label{lem:white}
  Suppose $L$ is a link, and $M$ is a (possibly disconnected) minimal surface
  spanning $L$.  Then for any point $p \in \R^3$ through which $n$
  sheets of~$M$ pass, the cone angle of~$L$ at $p$ is at least $2\pi n$.
\end{lemma}
\begin{proof}
  Let $S$ be the union of~$M$ and the exterior cone on $L$ from $p$.
  Consider the area ratio $\area(S\cap B_r(p))/(\pi r^2)$,
  where $B_r(p)$ is the ball of radius $r$ around~$p$ in~$\R^3$.
  As $r\to 0$, the
  area ratio approaches~$n$, the number of sheets of~$M$ passing through~$p$;
  as $r\to\infty$, the ratio approaches the density of the cone on~$L$ from~$p$,
  which is the cone angle divided by $2\pi$.  White has
  shown that the monotonicity formula for minimal surfaces continues
  to hold for $S$ in this setting~\cite{White-hild}:
  the area ratio is an increasing function of~$r$.
  Comparing the limit values at $r=0$ and $r=\infty$ we see that
  the cone angle from~$p$ is at least $2\pi n$.
\qed\end{proof}

As an immediate corollary, we obtain:

\begin{corollary}\label{cor:4pi}
  If $L$ is a nontrivial link, then there is some point $p$ from which
  $L$ has cone angle at least~$4\pi$.
\end{corollary}
\begin{proof}
  By the solution to the classical Plateau problem, each component
  of~$L$ bounds some minimal disk.  Let $M$ be the union of
  these disks.  Since $L$ is nontrivially linked, $M$ is not embedded: it must
  have a self-intersection point~$p$.  By the lemma, the cone angle at~$p$
  is at least~$4\pi$.
\qed\end{proof}

Note that, by Gauss--Bonnet, the cone angle of any cone over $K$ equals
the total geodesic curvature of~$K$ in the cone, which is clearly
bounded by the total curvature of~$K$ in space.  Therefore,
Corollary~\ref{cor:4pi} gives a new proof of the F\'ary--Milnor
theorem~\cite{fary,milnor}: any nontrivial link has total
curvature at least~$4\pi$.  (Compare~\cite[Cor.~2.2]{EWW}.)
This observation also shows
that the bound in Corollary~\ref{cor:4pi} cannot be improved,
since there exist knots with total curvature $4\pi+\epsilon$.

In fact any two-bridge knot can be built with total curvature
(and maximum cone angle) $4\pi+\epsilon$.
But we expect that for many knots of higher bridge number,
the maximum cone angle will necessarily be $6\pi$ or higher. 
For more information on these issues, see our paper~\cite{CKKS}
with Greg Kuperberg, where we give two
alternate proofs of Corollary~\ref{cor:4pi} in terms of the
\emph{second hull} of a link.

To apply the length estimate from Lemma~\ref{lem:pushoff}, we need a
stronger version for thick knots: If $K$ has thickness $\tau$, we must
show that the cone point of angle $4\pi$ can be chosen outside the
tube of radius $\tau$ surrounding $K$.

\begin{proposition}\label{prop:disk-outside-s}
  Let $K$ be a nontrivial knot, and let $\Tube$ be any embedded (closed) solid
  torus with core curve $K$. Any smooth disk $D$
  spanning $K$ must have self-intersections outside $\Tube$.
\end{proposition}
\begin{proof}
  Replacing $\Tube$ with a slightly bigger smooth solid torus if neccesary,
  we may assume that $D$ is transverse to the boundary torus
  $\bdy\Tube$ of~$\Tube$.  The intersection~$D\cap \bdy\Tube$ is then
  a union of closed curves.  If there is a self-intersection, we are
  done.  Otherwise,~$D \cap \bdy\Tube$ is a disjoint union of simple
  closed curves, homologous within~$\Tube$ to the core curve $K$ (via
  the surface $D\cap T$). Hence, within $\bdy\Tube$, its homology
  class $\alpha$ is the latitude plus some multiple of the meridian.
  Considering the possible arrangements of simple closed curves in the
  torus $\bdy\Tube$, we see that each intersection curve is homologous
  to zero or to $\pm\alpha$.

  Our strategy will be to first eliminate the trivial intersection curves,
  by surgery on $D$, starting with curves that are innermost on
  $\bdy\Tube$.  Then, we will find an essential intersection curve which is
  innermost on $D$: it is isotopic to $K$ and bounds a subdisk of~$D$
  outside $\Tube$, which must have self-intersections.

  To do the surgery, suppose $\gamma$ is an innermost intersection
  curve homologous to zero in $\bdy\Tube$.  It bounds a disk $A$
  within $\bdy\Tube$ and a disk $B$ within $D$.  Since $\gamma$ is an
  innermost curve on $\bdy\Tube$, $A \cap D$ is empty; therefore we
  may replace $B$ with $A$ without introducing any new
  self-intersections of~$D$.  Push~$A$ slightly off $\bdy\Tube$ to
  simplify the intersection.  Repeating this process a finite number
  of times, we can eliminate all trivial curves in $D \cap \bdy\Tube$.

  The remaining intersection curves are each homologous to $\pm\alpha$
  on $\bdy\Tube$ and thus isotopic to~$K$ within $\Tube$.  These do not
  bound disks on $\bdy\Tube$, but do on~$D$.  Some such curve $K'$ must be
  innermost on~$D$, bounding an open subdisk~$D'$.  Since~$K'$ is nontrivial
  in~$\Tube$, and $D' \cap \bdy\Tube$ is empty, the subdisk $D'$ must lie
  outside $\Tube$.  Because~$K'$ is knotted,~$D'$ must have
  self-intersections, clearly outside $\Tube$.
  Since we introduced no new self-intersections, these are self-intersections
  of $D$ as well.
\qed\end{proof}

We can now complete the proof of the main theorem of this section.

\begin{theorem}\label{thm:fourpi}
  If $K$ is a nontrivial knot then there is a point~$p$, outside the
  thick tube around~$K$, from which~$K$ has cone angle at
  least~$4\pi$.
\end{theorem}
\begin{proof}
  Span $K$ with a minimal disk~$D$, and
  let $T_n$ be a sequence of closed tubes around $K$,
  of increasing radius $\tau_n\to\tau(K)$.
  Applying Proposition~\ref{prop:disk-outside-s},
  $D$ must necessarily have a self-intersection point $p_n$ outside $T_n$.
  Using Lemma~\ref{lem:white}, the cone angle at~$p_n$ is at least $4\pi$.
  Now, cone angle is a continuous function on $\R^3$, approaching zero
  at infinity.  So the $p_n$ have a subsequence converging to some $p\in\R^3$,
  outside all the $T_n$ and thus outside the thick tube around $K$,
  where the cone angle is still at least $4\pi$.
\qed\end{proof}

It is interesting to compare the cones of cone angle $4\pi$
constructed by Theorem~\ref{thm:fourpi} with those of cone angle $2\pi$
constructed by Lemma~\ref{lem:twopi}; see Figure~\ref{fig:3-1-4pi}.

\figtwox{3-1-4pi}{1.1}{1.8}{A $4\pi$ cone on a trefoil knot.}
{Two views of a cone, whose cone angle is precisely $4\pi$,
on the symmetric trefoil knot from Figure~\ref{fig:3-1-2pi}.
Computational data shows this is close
to the maximum possible cone angle for this trefoil.}

\section{Parallel Overcrossing Number Bounds for Knots}\label{sec:newbounds}

We are now in a position to get a better lower bound for the ropelength
of any nontrival knot.

\begin{theorem}\label{thm:newbounds}
  For any nontrivial knot $K$ of unit thickness,
  \[
    \Len(K) \geq 4\pi +2\pi \sqrt{\pov(K)} \geq 2\pi\big(2+\sqrt2\big).
  \]
\end{theorem}
\begin{proof}
  Let $\Tube$ be the thick tube (the unit-radius solid torus)
  around~$K$, and let $V$ be the $C^1$ unit vectorfield inside~$\Tube$ as in
  the proof of Theorem~\ref{thm:link}.  Using
  Theorem~\ref{thm:fourpi}, we construct a cone surface $S$ of cone
  angle~$4\pi$ from a point $p$ outside $T$.

  Let $S'$ be the cone defined by deleting a unit neighborhood
  of $\bdy S$ in the intrinsic geometry of $S$. Take any $q \in K$
  farthest from the cone point $p$.  The intersection of~$S$ with the
  unit normal disk~$D$ to~$K$ at~$q$ consists only of the unit line
  segment from $q$ towards $p$; thus $D$ is disjoint from $S'$.

  In general, the integral curves of $V$ do not close. However, we can
  define a natural map from $\Tube \setm D$ to the unit disk $D$ by
  flowing forward along these integral curves. This map is continuous
  and distance-decreasing.  Restricting it to $S'\cap T$ gives a
  distance-decreasing (and hence area-decreasing) map to $D$, which we
  will prove has unsigned degree at least~$\pov(K)$.

  \figside{inter-and-proj}{.32}{The projection from $S' \cap T$ to $D$.}
  {This trefoil knot, shown with its thick tube~$T$, is coned to the
  point $p$ to form the cone surface $S'$, as in the proof of
  Theorem~\ref{thm:newbounds}.  The disk $D$ is normal to the knot at
  the point furthest from $p$.  We follow two integral curves of $V$
  within $T\setm D$, through at least~$\pov(K)$ intersections with
  $S'$, until they end on $D$. Although we have drawn the curves as if
  they close after one trip around~$T$, this is not always the case.}{.2}{.7}

  Note that $K' := \bdy S'$ is isotopic to $K$ within $\Tube$, and thus
  $\pov(K) = \pov(K')$. Furthermore, each integral curve~$C$ of~$V$ in
  $\Tube \setm D$ can be closed by an arc within $D$ to a knot $C'$
  parallel to $K'$. In the projection of~$C'$ and $K'$ from the
  perspective of the cone point, $C'$ must overcross $K'$ at least
  $\pov(K')$ times. Each of these crossings represents an intersection
  of~$C'$ with $S'$. Further, each of these intersections is an
  intersection of~$C$ with $S'$, since the portion of~$C'$ not in $C$
  is contained within the disk $D$. This proves that our
  area-decreasing map from $S' \cap \Tube$ to $D$ has unsigned degree at least
  $\pov(K)$. (An example of this map is shown in
  Figure~\ref{fig:inter-and-proj}.) Since~$\Area(D) = \pi$ it follows
  that 
  \[ 
  \Area(S') \geq \Area(S' \cap \Tube) \geq \pi \pov(K).  
  \]

  The isoperimetric inequality in a $4\pi$ cone is affected by the
  negative curvature of the cone point.  However, the length $\ell$ required
  to surround a fixed area on $S'$ is certainly no less than that required in
  the Euclidean plane:
  \[\label{eqn:cone-k'}
    \ell \geq 2\pi \sqrt{\pov(K)}.
  \] 
  Since each point on $K'$ is at unit distance from $K$, we
  know $S'$ is surrounded by a unit-width neighborhood inside
  $S$. Applying Lemma~\ref{lem:pushoff} we see that 
 \[ 
   \Len(K) \geq 4\pi +2\pi \sqrt{\pov(K)},
  \] 
  which by Proposition~\ref{prop:bridge} is at least $2\pi(2+\sqrt2)$.
\qed\end{proof}

\section{Asymptotic Crossing Number Bounds for Knots and Links}
\label{sec:acbounds}

The proof of Theorem~\ref{thm:newbounds} depends on the fact that $K$
is a single knot: for a link~$L$, there would be no guarantee that we
could choose spanning disks~$D$ for the tubes around the components of~$L$
which were all disjoint
from the truncated cone surface.  Thus, we would be unable to close
the integral curves of~$V$ without (potentially) losing crossings in
the process.

We can overcome these problems by using the notion of asymptotic crossing
number.  The essential idea of the proof is that (after a small
deformation of~$K$) the integral curves of~$V$ will close after some
number of trips around~$K$. We will then be able to complete the proof
as above, taking into account the complications caused by traveling
several times around~$K$.

For a link $L$ of~$k$ components, $K_1, \dots, K_k$, Freedman and He
\cite{Freedman-He} define a relative asymptotic crossing number
\[
\ac(K_i,L) := \inf_{pK_i,qL} \frac{\Over(pK_i,qL)}{|pq|},
\]
where the infimum is taken over all degree-$p$ satellites $pK_i$ of~$K_i$ and
all degree-$q$ satellites $qL$ of $L$.  It is easy to see that, for
each~$i$,
\[
\ac(K_i,L) \geq \sum_{j \neq i} |\Lk(K_i,K_j)|.
\]

Freedman and He also give lower bounds for this asymptotic crossing
number in terms of genus, or more precisely the Thurston norm.  To
understand these, let $T$ be a tubular neighborhood of~$L$.  Then
$H_1(\bdy T)$ has a canonical basis consisting of latitudes $l_i$ and
meridians $m_i$.  Here, the latitudes span the kernel of the map
$H_1(\bdy T) \to H_1(\R^3 \setm T)$ induced by inclusion, while the
meridians span the kernel of $H_1(\bdy T) \to H_1(T)$.

The boundary map $H_2(\R^3 \setm T, \bdy T) \to H_1(\bdy T)$ is an
injection; its image is spanned by the classes
\[
\alpha_i := l_i + \sum_{j\neq i} \Lk(K_i,K_j) m_j.
\]
We now define 
\[
\chi_-(K_i,L) := \min_S |S|_T,
\]
where the minimum is taken over all embedded surfaces $S$
representing the (unique) preimage of $\alpha_i$ in
$H_2(\R^3 \setm T, \bdy T)$, and $|S|_T$ is the Thurston norm
of the surface $S$.  That is,
$$|S|_T := \sum_{S_k} -\chi(S_k),$$
where the sum is taken over all components $S_k$ of $S$ which
are not disks or spheres, and $\chi$ is the Euler characteristic.
With this definition, Freedman and He prove~\cite[Thm.~4.1]{Freedman-He}:
\begin{proposition} \label{prop:fh}
If $K$ is a component of a link $L$,
\[
\ac(K,L) \geq \chi_{-}(K,L).
\]
In particular, 
$\ac(K,L) \geq 2 \genus(K) - 1$.
\qed\end{proposition}

Our interest in the asymptotic crossing number comes from the following bounds:
\begin{theorem}
\label{thm:acbounds}
Suppose $K$ is one component of a link~$L$ of unit thickness. Then 
\[ \Len(K) \geq 2\pi +2\pi\sqrt{\ac(K,L)}. \]
If $K$ is nontrivially knotted, this can be improved to 
\[ \Len(K) \geq 4\pi +2\pi\sqrt{\ac(K,L)}. \]
\end{theorem}
\begin{proof}
As before, we use Lemma~\ref{lem:twopi} or Theorem~\ref{thm:fourpi} to
construct a cone surface $S$ of cone angle $2\pi$ or $4\pi$.
We let $S'$ be the complement of a unit neighborhood of $\bdy S$,
and set $K':=\bdy S'$, isotopic to~$K$.

Our goal is to bound the area of~$S'$ below. As before, take the
collection~$\Tube$ of embedded tubes surrounding the components
of~$L$, and let $V$ be the $C^1$ unit vectorfield normal to the normal disks
of~$\Tube$. Fix some component~$J$ of~$L$ (where $J$ may be the same
as $K$), and any normal disk $D$ of the embedded tube~$\Tube_J$ around~$J$. The
flow of~$V$ once around the tube defines a map from $D$ to $D$. The
geometry of~$V$ implies that this map is an isometry, and hence this
map is a rigid rotation by some angle $\theta_J$. Our first claim is
that we can make a $C^1$-small perturbation of~$J$ which ensures that
$\theta_J$ is a rational multiple of~$2\pi$.

Fix a particular integral curve of~$V$. Following this integral curve
once around~$J$ defines a framing (or normal field) on $J$ which fails
to close by angle $\theta_J$. If we define the {\em twist} of a
framing $W$ on a curve $J$ by
\[
\Tw(W) = \frac{1}{2\pi} \int \frac{dW}{ds} \cross W \cdot \ds,
\]
it is easy to show that this framing has zero twist.  We can close this
framing by adding twist $-\theta_J/2\pi$, defining a framing $W$ on $J$.
If we let $\Wr(J)$ be the \emph{writhe} of~$J$, then the
C{\u{a}}lug{\u{a}}reanu--White formula \cite{cal1,cal2,cal3,white}
tells us that $\Lk(J,J') = \Wr(J) - \theta_J/2\pi$,
where $J'$ is a normal pushoff of~$J$ along~$W$.  Since the linking
number $\Lk(J,J')$ is an integer, this means that
$\theta_J$ is a rational multiple of $2\pi$
if and only if $\Wr(J)$ is rational.
But we can alter the writhe of~$J$ to be rational with a $C^1$-small
perturbation of~$J$ (see \cite{Fuller,Benham} for details), proving
the claim.

So we may assume that, for each component $J$ of~$L$,
$\theta_J$ is a rational multiple $2\pi p_J/q_J$ of $2\pi$.
Now let $q$ be the least common multiple of the (finitely many) $q_J$.
We will now define a distance- and area-decreasing map of
unsigned degree at least $q\ac(K,L)$ from the intersection of $T$ and
the cone surface $S'$ to a sector of the unit disk of angle $2\pi/q$.

Any integral curve of $V$ must close after $q_J$ trips around~$J$.
Thus, the link~$J^q$ defined by following the integral curves
through $q/q_J$ points spaced at angle $2\pi/q$ around a normal disk
to~$J$ is a degree-$q$ satellite of~$J$.
Further, if we divide a normal disk to~$J$ into sectors of angle $2\pi/q$,
then $J^q$ intersects each sector once.

We can now define a distance-decreasing map from $S' \cap T_J$ to the
sector by projecting along the integral curves of~$V$. Letting $L^q$
be the union of all the integral curves $J^q$, and identifying the
image sectors on each disk gives a map from $S' \cap T_L$ to the
sector. By the definition of $\ac(K,L)$,
\[
  \Over(L^q,K') = \Over(K,L^q) \geq q\ac(K,L),
\]
so $L^q$ overcrosses $K'$ at least $q\ac(K,L)$ times.  Thus we have at least
$q\ac(K,L)$ intersections between~$L^q$ and~$S'$, as in the proof of
Theorem~\ref{thm:newbounds}. Since the sector has area $\pi/q$, this
proves that the cone $S'$ has area at least~$\pi\ac(K,L)$, and thus
perimeter at least $2\pi\sqrt{\ac(K,L)}$. The theorem
then follows from Lemma~\ref{lem:pushoff} as usual.
\qed\end{proof}

Combining this theorem with Proposition~\ref{prop:fh} yields:
\begin{corollary}\label{cor:chibounds}
For any nontrivial knot~$K$ of unit thickness,
\[
\Len(K) \geq 2\pi\Big(2 + \sqrt{2 \genus(K) - 1}\Big).
\]
For any component~$K$ of a link~$L$ of unit thickness,
\[ \hfill \Len(K) \geq 2\pi\big(1 + \sqrt{\chi_-(K,L)}\big), \]
where $\chi_-$ is the minimal Thurston norm as above.
\end{corollary}

As we observed earlier, $\ac(K_i,L)$ is at least the sum of the linking numbers
of the~$K_j$ with~$K_i$, so Theorem~\ref{thm:acbounds} subsumes
Theorem~\ref{thm:link}. Often, it gives more information. When the
linking numbers of all $K_i$ and $K_j$ vanish, the minimal Thurston
norm $\chi_{-}(K_i,L)$ has a particularly simple interpretation: it is
the least genus of any embedded surface spanning $K_i$ and
avoiding~$L$.  For the Whitehead link and Borromean rings, this
invariant equals one, and so these bounds do not provide an
improvement over the simple-minded bound of Theorem~\ref{thm:peribound}.

To find an example where Corollary~\ref{cor:chibounds} is an improvement,
we need to be able to compute the Thurston norm.
McMullen has shown~\cite{McMullen} that the Thurston norm is bounded
below by the Alexander norm, which is easily computed from the multivariable
Alexander polynomial.
One example he suggests is a $(2,2n)$--torus link with two components.
If we replace one component $K$ by its Whitehead double, then in the new link,
the other component has Alexander norm $2n-1$.  Since it is clearly spanned
by a disk with $2n$ punctures (or a genus~$n$ surface) avoiding~$K$,
the Thurston norm is also $2n-1$.  Figure~\ref{fig:thurston} (left) shows
the case $n=3$, where the Alexander polynomial is $(1+x+x^2)^2(1-x)(1-y)$.

On the other hand, if $K$ is either component of the three-fold link $L$
on the right in Figure~\ref{fig:thurston}, we can span $K$ with a genus-two
surface, so we expect that $\chi_{-}(K,L) = 3$,
which would also improve our ropelength estimate.  However, it seems hard
to compute the Thurston norm in this case.  The Alexander norm in
this case is zero, and even the more refined bounds of Harvey~\cite{Harvey}
do not show the Thurston norm is any greater.


\figtwo{thurston}{1.1}{Links with possibly large Thurston norm.}
{At the left we see the result of replacing one component
of a $(2,6)$--torus link by its Whitehead double.  In this link $L$,
the other component has Alexander norm, and hence also Thurston norm,
equal to $5$.  Thus Corollary~\ref{cor:chibounds} shows the total ropelength
of $L$ is at least $2\pi(3+\sqrt5)$.
For the three-fold link at the right, which is a bangle sum of three
square-knot tangles, we expect the Thurston norm to be $3$
(which would give ropelength at least $4\pi(1+\sqrt3)$),
but we have not found a way to prove it is not less.}

\section{Asymptotic Growth of Ropelength}\label{sec:asymptotic}

All of our lower bounds for ropelength have been asymptotically proportional
to the square root of the number of components, linking number,
parallel crossing number, or asymptotic crossing number.
While our methods here provide the best known results for fairly small links,
other lower bounds grow like the $\tfrac34$ power of these complexity
measures.  These are of course better for larger links,
as described in our paper
{\em Tight Knot Values Deviate from Linear Relation}~\cite{CKS-nature}.
In particular, for a link type~$\Ltype$ with crossing number~$n$,
the ropelength is at least $\big(\tfrac{4\pi}{11}n\big)^{3/4}$,
where the constant comes from~\cite{BuckSimon}.
In~\cite{CKS-nature} we gave examples
(namely the $(k,k-1)$--torus knots and the $k$-component Hopf links,
which consist of $k$ circles from a common Hopf fibration of~$\S^3$)
in which ropelength grows exactly as the $\tfrac34$ power of crossing number.

Our Theorem~\ref{thm:peribound} proves that for the
the simple chains (Figure~\ref{fig:tight-chain}),
ropelength must grow linearly in crossing number~$n$.
We do not know of any examples exhibiting superlinear growth,
but we suspect they might exist, as described below.

To investigate this problem, consider representing a link type~$\Ltype$ with
unit edges in the standard cubic lattice~$\Z^3$.  The minimum number of
edges required is called the lattice number~$k$ of~$\Ltype$.  We claim
this is within a constant factor of the ropelength~$\ell$ of a tight
configuration of~$\Ltype$.
Indeed, given a lattice representation with $k$ edges,
we can easily round off the corners with quarter-circles of radius~$\half$
to create a $C^{1,1}$ curve with length less than $k$ and thickness~$\half$,
which thus has ropelength~$\ell$ at most~$2k$.
Conversely, it is clear that any thick knot of ropelength~$\ell$
has an isotopic inscribed polygon with $\Order(\ell)$ edges and
bounded angles; this can then be replaced by an isotopic lattice
knot on a sufficiently small scaled copy of~$\Z^3$.  We omit our detailed
argument along the lines, showing $k\le 94\ell$,
since Diao \emph{et al.}~\cite{DEJvR-ln} have
recently obtained the better bound $k\le12\ell$.

The lattice embedding problem for links is similar to the
VLSI layout problem~\cite{leighton,leighton-lb},
where a graph whose vertex degrees are at most~$4$
must be embedded in two layers of a cubic lattice.
It is known~\cite{BhaLei} that any $n$-vertex planar graph can be
embedded in VLSI layout area $\Order\big(n(\log n)^2\big)$.
Examples of planar $n$-vertex graphs requiring layout area
at least $n\log n$ are given by the so-called trees of meshes.
We can construct $n$-crossing links analogous to these trees of meshes,
and we expect that they have lattice number at least~$n\log n$,
but it seems hard to prove this.
Perhaps the VLSI methods can also be used to show that lattice number
(or equivalently, ropelength) is at most $\Order\big(n(\log n)^2\big)$.

Here we will give a simple proof that the ropelength of an $n$-crossing
link is at most~$24n^2$, by constructing a lattice embedding of
length less than~$12n^2$.  This follows from the
theorem of Schnyder~\cite{dFPP,Schnyder}
which says that an $n$-vertex planar graph
can be embedded with straight edges connecting
vertices which lie on an $(n-1)\times (n-1)$ square grid.  We double this
size, to allow each knot crossing to be built on a $2\times2\times2$
array of vertices.  For an $n$-crossing link diagram, there are
$2n$ edges, and we use $2n$ separate levels for these edges.
Thus we embed the link in a $(2n-2)\times(2n-2)\times(2n+2)$
piece of the cubic lattice.  Each edge has length less than $6n$,
giving total lattice number less than $12n^2$.

Note that Johnston has recently given an independent
proof~\cite{heather} that an $n$-crossing
knot can be embedded in the cubic grid with length $\Order(n^2)$.
Although her constant is worse than our $12$, her embedding is
(like a VLSI layout) contained in just two layers of the cubic lattice.
It is tempting to think that an~$\Order(n^2)$ bound
on ropelength could be deduced from the Dowker code for a knot,
and in fact such a claim appeared in~\cite{Buck-nature}.
But we do not see any way to make such an argument work.

The following theorem summarizes the results of this section:
\begin{theorem}
Let $\Ltype$ be a link type with minimum crossing number~$n$,
lattice number~$k$, and minimum ropelength $\ell$.  Then
\[ \left(\frac{4\pi}{11} n\right)^{3/4}
    \,\le\, \ell \,\le\, 2k \,\le\, 24 n^2. \qquad\qed\]
\end{theorem}

\section{Further Directions}\label{sec:conclusion}

Having concluded our results, we now turn to some open problems and conjectures.

The many examples of tight links constructed in
Section~\ref{sec:peribound} show that the existence and regularity
results of Section~\ref{sec:exist} are in some sense optimal: we know that
ropelength minimizers always exist, we cannot expect a ropelength
minimizer to have global regularity better than $C^{1,1}$,
and we have seen that there exist
continuous families of ropelength minimizers with different shapes.
Although we know that each ropelength minimizer has well-defined
curvature almost everywhere (since it is $C^{1,1}$) it would be
interesting to determine the structure of the singular set where
the curve is not $C^2$. We expect this singular set is finite,
and in fact:

\begin{conjecture} 
Ropelength minimizers are piecewise analytic.
\end{conjecture}

The $\Peri{n}$ bound for the ropelength of links in
Theorem~\ref{thm:peribound} is sharp, and so cannot be improved. But
there is a certain amount of slack in our other ropelength
estimates. The parallel crossing number and asymptotic crossing number
bounds of Section~\ref{sec:newbounds} and Section~\ref{sec:acbounds}
could be immediately improved by showing:

\begin{conjecture}
If $L$ is any knot or link, $\ac(L) = \pov(L) = \MCN(L)$.
\end{conjecture}

For a nontrivial knot, this would increase our best estimate to $4\pi +
2\pi\sqrt{3} \approx 23.45$, a little better than our current estimate
of $4\pi + 2\pi\sqrt{2} \approx 21.45$ (but not good enough to decide
whether a knot can be tied in one foot of one-inch rope). A more
serious improvement would come from proving:

\begin{conjecture}\label{conj:disj}
The intersection of the tube around a knot of unit thickness with some
$4\pi$ cone on the knot contains $\pov(K)$ disjoint unit disks
avoiding the cone point.
\end{conjecture}

Note that the proof of Theorem~\ref{thm:newbounds} shows only that this
intersection has the \emph{area} of $\pov(K)$ disks. This conjecture
would improve the ropelength estimate for a nontrivial knot to about
$30.51$, accounting for $93\%$ of Pieranski's numerically computed
value of $32.66$ for the ropelength of the trefoil~\cite{Pieranski}.
We can see the tightness of this proposed estimate in
Figure~\ref{fig:ideal-4pi-inter}.%
\figh{ideal-4pi-inter}{1.5}{The intersection of the tube on
Pieranski's ideal trefoil with one-third of a $4\pi$ cone on that
trefoil, developed into the plane.}
{Pieranski's numerically computed tight trefoil $K$ has three-fold symmetry,
and there is a $4\pi$ cone point~$p$ on the symmetry axis.  The cone from
$p$ also has three-fold symmetry, and a fundamental sector
develops into a $4\pi/3$ wedge in the plane, around the point~$p$.
Here we show the development of that sector.  The shaded regions
are the intersection of the cone with the thick tube around~$K$.
These include a strip (of width at least~$1$) inwards from the
boundary $K$ of the cone, together with a disk around
the unique point $q$ where $K$ cuts this sector of the cone.
Our Conjecture~\ref{conj:disj} estimates the area of the cone from below
by the area of a unit disk around $q$ plus a unit-width strip around $K$.
The figure shows that the actual shaded disk and strip are not much
bigger than this, and that they almost fill the sector.}

Very recently, Diao has announced~\cite{Diao} a proof that the length
of any unit-thickness knot $K$ satisfies
$$16\pi\MCN(K) \le \Len(K)\big(\Len(K)-17.334\big).$$
This improves our bounds in many cases.  He also finds that
the ropelength of a trefoil knot is greater than $24$.

Our best current bound for the ropelength of the Borromean rings is
$12\pi \approx 37.70$, from Theorem~\ref{thm:peribound}.  Proving only
the conjecture that $\ac(L) = \MCN(L)$ would give us a fairly sharp
bound on the total ropelength: If each component has asymptotic
crossing number~$2$, Theorem~\ref{thm:acbounds} tells us that
$6\pi(1+\sqrt{2}) \approx 45.51$
is a bound for ropelength. This bound would account
for at least $78\%$ of the optimal ropelength, since we can exhibit a
configuration with ropelength about $58.05$, built from three congruent
planar curves, as in Figure~\ref{fig:brings-two}.
\figthree{brings-two}{1.4}{Borromean rings.}
{This configuration of the Borromean rings has ropelength about $58.05$.
It is built from three congruent piecewise-circular plane curves,
in perpendicular planes.  Each one consists of arcs
from four circles of radius~$2$ centered at the vertices of a rhombus
of side~$4$, whose major diagonal is $4$ units
longer than its minor diagonal.}

Although it is hard to see how to improve the ropelength
of this configuration of the Borromean rings, it is not tight.
In work in progress with Joe Fu, we define a notion of criticality
for ropelength, and show that this configuration
is not even ropelength-critical.

Finally, we observe that our cone surface methods seem useful in many
areas outside the estimation of ropelength. For example,
Lemma~\ref{lem:twopi} provides the key to a new proof an unfolding
theorem for space curvess:

\begin{proposition}
\label{prop:unfolding}
For any space curve $K : S^1 \to \R^3$, parametrized by arclength,
there is a plane curve $K'$ of the same length, also parametrized by
arclength, so that for every $\theta$, $\phi$ in $S^1$,
\[
|K(\theta) - K(\phi)| \leq |K'(\theta) - K'(\phi)|.
\]
\end{proposition}
\begin{proof}
By Lemma~\ref{lem:twopi}, there exists some cone point $p$ for which the cone 
of $K$ to $p$ has cone angle $2\pi$. Unrolling the cone on the plane,
an isometry, constructs a plane curve $K'$ of the same arclength. Further, each
chord length of $K'$ is a distance measured
in the instrinsic geometry of the cone,
which is at least the corresponding distance in $\R^3$.
\qed\end{proof}

This result was proved by Reshetnyak~\cite{Resh1,Resh2} in a more general
setting: a curve in a metric space of curvature bounded above (in the sense
of Alexandrov) has an unfolding into the model two-dimensional space of
constant curvature.
The version for curves in Euclidean space was also proved
independently by Sallee~\cite{Sal}.
(In~\cite{KS-disto}, not knowing of this earlier work, we stated
the result as Janse van Rensburg's unfolding conjecture.)

The unfoldings of Reshetnyak and Sallee are always convex curves in the plane.
Our cone surface method, given in the proof of Proposition~\ref{prop:unfolding},
produces an unfolding that need not be convex, as shown in
Figure~\ref{fig:unfolding}. Ghomi and Howard have recently extended
our argument to prove stronger results about unfoldings~\cite{GH}.

\figh{unfolding}{1.6}{A trefoil knot and its planar unfolding.}  {On
the left, we see a trefoil knot, with an example chord. On the right
is a cone-surface unfolding of that knot, with the corresponding
chord.  The pictures are to the same scale, and the chord on the right
is clearly longer.}

\begin{acknowledgement}
We are grateful for helpful and productive conversations with
many of our colleagues, including Uwe Abresch, Colin Adams,
Ralph Alexander, Stephanie Alexander, Therese Biedl, Dick Bishop,
Joe Fu, Mohammad Ghomi, Shelly Harvey, Zheng-Xu He, Curt McMullen,
Peter Norman, Saul Schleimer and Warren Smith.
We would especially like to thank Piotr Pieranski for providing the data
for his computed tight trefoil,
Brian White and Mike Gage for bringing Lemma~\ref{lem:white}
and its history to our attention,
and the anonymous referee for many detailed and helpful suggestions.
Our figures were produced with Geomview and Freehand.
This work has been supported by the National Science Foundation
through grants DMS-96-26804 (to the GANG lab at UMass),
DMS-97-04949 and DMS-00-76085 (to Kusner),
DMS-97-27859 and DMS-00-71520 (to Sullivan),
and through a Postdoctoral Research Fellowship DMS-99-02397 (to Cantarella).
\end{acknowledgement}

\newcommand{\etalchar}[1]{$^{#1}$}

\end{document}